\documentclass[12pt] {article}
\usepackage{amsmath,amsthm}
\usepackage{amssymb,latexsym}

\usepackage{empheq}
\usepackage{amscd}
\usepackage{chngcntr}
\usepackage{enumitem}
\newlist{paragraphlist}{enumerate}{1}

\setlist[paragraphlist,1]{leftmargin=*,label={\bfseries \arabic*}}

\counterwithin{paragraphlisti}{subsubsection}

\title{Kotrbat\'y's theorem on valuations and geometric inequalities for convex bodies.}
\date{}
\author{Semyon Alesker \footnote{Partially supported by ISF grant 865/16 and the  US - Israel BSF grant 2018115.}
\\  { \normalsize Department of Mathematics, Tel Aviv University, Ramat Aviv}
\\  { \normalsize 69978 Tel Aviv, Israel }
\\ {\normalsize e-mail: semyon@tauex.tau.ac.il }}

\def\RR{\mathbb{R}}
\def\CC{\mathbb{C}}

\def\eps{\varepsilon}
\def\alp{\alpha}

\def\lam{\lambda}

\def\to{\longrightarrow}

\def\qed { Q.E.D. }

\def\inj{\hookrightarrow}

\swapnumbers
\newtheorem{theorem}{Theorem}[section]
\newtheorem{corollary}[theorem]{Corollary}
\newtheorem{lemma}[theorem]{Lemma}

\theoremstyle{definition}

\newtheorem{example}[theorem]{Example}
\newtheorem{definition}[theorem]{Definition}
\newtheorem{remark}[theorem]{Remark}
\theoremstyle{conjecture}
\newtheorem{conjecture}[theorem]{Conjecture}
\theoremstyle{principle}

 \def\ck{{\cal K}}

\def\pt{\partial}

\numberwithin{equation}{section}
\begin{document}
\maketitle

\begin{abstract}
Very recently J. Kotrbat\'y [13] has proven general inequalities for
translation invariant smooth valuations formally analogous to the Hodge-Riemann bilinear relations in the K\"ahler geometry. The goal of this
note is to prove  several inequalities for mixed volumes of convex bodies using Kotrbat\'y's theorem as the main tool.
\end{abstract}

\section{Introduction.}\label{S:introduction}
\begin{paragraphlist}
\item Very recently J. Kotrbat\'y [13] has proven general inequalities for translation invariant smooth valuations analogous to the Hodge-Riemann bilinear relations from the K\"ahler geometry, see Theorem 2.14 below.
The goal of this note is prove several new inequalities for mixed volumes of convex bodies using his theorem as the main tool.

To formulate the main results let us introduce a notation. Let $\RR^n$ denote the standard Euclidean space of dimension $n$. Let us denote by  $\iota_1,\iota_2,\Delta\colon\RR^n\to \RR^{2n}=\RR^n\times \RR^n$ the maps
$$\iota_1(x)=(x,0),\, \iota_2(x)=(0,x),\, \Delta(x)=(x,x).$$
We refer to the book \cite{schneider-2ed} for the notions of mixed volumes and mixed area measures.
The main result of this paper is
\begin{theorem}\label{T:main-intro}
Let $n\geq 2$. Let $A_1,\dots,A_{n-1}\subset\RR^n$ be convex compact subsets.
Let $B\subset \RR^n$ be the unit Euclidean ball.

(1) (a) Then one has inequality for mixed volumes in $\RR^{2n}=\RR^n\times\RR^n$
\begin{eqnarray*}\label{E:odd-ineq}
V(\iota_1A_1,\dots,\iota_1A_{n-1};\iota_2A_1,\dots,\iota_2A_{n-1};\Delta(B)[2])\geq\\ V(\iota_1A_1,\dots,\iota_1A_{n-1};-\iota_2A_1,\dots,-\iota_2A_{n-1};\Delta(B)[2]).
\end{eqnarray*}
(b) In the special case when $A_1=\dots=A_{n-1}=:A$ with $A$ having non-empty interior and smooth boundary with positive Gauss curvature the equality in the above inequality
is achieved if and only if $A$ has a center of symmetry.

(2) (a) Furthermore
\begin{eqnarray*}
V(\iota_1A_1,\dots,\iota_1A_{n-1};\iota_2A_1,\dots,\iota_2A_{n-1};\Delta(B)[2])+\\
+V(\iota_1A_1,\dots,\iota_1A_{n-1};-\iota_2A_1,\dots,-\iota_2A_{n-1};\Delta(B)[2])\leq \gamma_n V(A_1,\dots,A_{n-1},B)^2,
\end{eqnarray*}
where $\gamma_n$ is a constant depending on $n$ only and uniquely characterized by the property that if $A_1=\dots=A_{n-1}=B$ then in the above inequality there is an equality.

(b) In the special case when $A_1=\dots=A_{n-1}=:A$ with $A$ being centrally symmetric, having non-empty interior and smooth boundary with positive Gauss curvature the equality in the above inequality
is achieved if and only if $A$ is a Euclidean ball of arbitrary radius.

(3) Let $n\geq 4$. Let $K_1,K_2,L_1,L_2\subset \RR^n$ be centrally symmetric convex compact sets with smooth boundary with positive Gauss curvature.
Let us assume the following equality of the mixed area measures:
\begin{eqnarray}\label{E:mixed-meas-cond-in}
S(K_1,K_2,B[n-3],\cdot)=S(L_1,L_2,B[n-3],\cdot).
\end{eqnarray}
Then one has
\begin{eqnarray}\label{E:ineq-2-homog}
V(K_1[2],K_2[2],B[n-4])+V(L_1[2],L_2[2],B[n-4])\geq 2V(K_1,K_2,L_1,L_2,B[n-4]).
\end{eqnarray}
\end{theorem}

\begin{remark}\label{remark-main-thm}
(1) It is expected that in part (3) the assumption that all bodies $K_i,L_i$ are centrally symmetric is unnecessary. That would follow from Kotrbat\'y's conjecture \ref{kotrbaty-conj}.

(2) If in part (3) of the theorem one takes $K_1=K_2=:K$ and $L_1=L_2=:L$ then the inequality (\ref{E:ineq-2-homog}) is weaker than what is actually known: in this case the equality is achieved.
Indeed the Alexandrov-Fenchel-Jessen theorem (see \cite{schneider-2ed}, Corollary 8.1.4) says that the assumption (\ref{E:mixed-meas-cond-in}) implies that $K$ and $L$ have to be translates of each other.

(3) In part (2)(b) the assumption of central symmetry of $A$ cannot be omitted at least in the case of $n=2$. One can show that for $n=2$ the inequality (2)(a) is equivalent to the standard isoperimetric inequality for the body $A-A$:\footnote{This observation is a result of 
a correspondence with T. Wannerer.}
$$area(A-A)\leq \frac{length(\pt (A-A))^2}{4\pi}=\frac{length(\pt A)^2}{\pi}.$$
The equality here is achieved if and only if $A-A$ is a Euclidean disk of arbitrary radius.
\end{remark}

\item Valuations on convex sets are finitely additive functionals on convex compact subsets in $\RR^n$ (the  Definition \ref{D:valuation}). It is a classical object of convexity.
During the last 25 years there was a considerable progress in the theory.
There was a breakthrough in 1995 when the Klain-Schneider characterization of translations invariant simple continuous valuations was obtained \cite{klain}, \cite{schneider-simple}
which opened a way to further progress. In particular the McMullen's conjecture has been solved \cite{alesker-mcmullen}. Subsequent progress led to a discovery of a few new structures on valuations, to posing new questions, and to opening new
directions of the research. Some of these structures turned out to be useful in integral geometry, see e.g. \cite{bernig-fu-annals}, \cite{bernig-fu-solanes}.
 Examples of new structures relevant to this paper are the product on smooth valuations introduced by the author \cite{alesker-product} and
the convolution on them introduced by Bernig and Fu \cite{bernig-fu-convolution} (see Theorems \ref{T:product} and \ref{T;convolution} below).

The most basic  examples of
continuous translation invariant valuations are given by mixed volumes of convex sets. Mixed volume is a well studied classical object in convexity of its own right.
One of the most remarkable properties of the mixed volume is the Alexandrov-Fenchel inequality (see \cite{schneider-2ed}, Theorem 7.3.1). Kotrbat\'y \cite{kotrbaty-hr} has stated a conjecture,
which is a generalized version of a part of his Conjecture \ref{kotrbaty-conj} below, which implies the Alexandrov-Fenchel inequality. Furthermore he deduced there a special case of the latter inequality from
his Theorem \ref{T:kotrbaty} below.

The main result of this paper, Theorem \ref{T:main-intro}, shows that Kotrbat\'y's work \cite{kotrbaty-hr} can also be used to obtain new geometric inequalities using valuations theoretical methods.
Modulo some background from convexity in general and valuations theory in particular, Theorem \ref{T:main-intro} is rather straightforward consequence of the
Kotrbat\'y's Theorem \ref{T:kotrbaty} which is applied to certain specific valuations and the result is interpreted in terms of mixed volumes.

\item In Section \ref{S:background} we review some background from the valuations theory necessary to formulate the Kotrbat\'y's theorem.
Section \ref{S:pf-main} contains proofs of the main results. Theorem \ref{T:main-intro} summarizes Theorems \ref{T:odd-ineq}, \ref{T:even-ineq}, and \ref{T:mixed-area} in the main text.

\item {\bf Acknowledgements.} I thank J. Kotrbat\'y for useful remarks on the first version of the paper. I thank R. Schneider and T. Wannerer for useful correspondences.
\end{paragraphlist}
\section{Background.}\label{S:background}
\begin{paragraphlist}
\item In this section we summarize  basic definitions and results on valuations. For general background on convexity and, in particular, mixed volumes we refer to \cite{schneider-2ed}.

For an $n$-dimensional real vector space $W$ we denote by $\ck(W)$ the family of all convex compact non-empty subsets of $W$. Being equipped with the Hausdorff metric $\ck(W)$
is a locally compact space by the Blaschke selection theorem (\cite{schneider-2ed}, Theorem 1.8.7).

\item
\begin{definition}\label{D:valuation}
A {\itshape valuation} is a map $\phi\colon \ck(W)\to \CC$ which satisfies the following additivity property:
$$\phi(A\cup B)=\phi(A)+\phi(B)-\phi(A\cap B)$$
whenever $A,B,A\cup B\in\ck(W)$.
\end{definition}

\begin{definition}\label{D:contin-val}
A valuation $\phi$ is called {\itshape continuous} if it is continuous in the Hausdorff metric.
\end{definition}

\item
\begin{definition}\label{D:transl-invar}
A valuation $\phi$ is called {\itshape translation invariant} if
$$\phi(K+x)=\phi(K)$$
for any $K\in\ck(W),x\in W$.
\end{definition}

We will denote by $Val(W)$ the space of continuous translation invariant valuations on $\ck(W)$. Being equipped with the topology of uniform convergence on compact subsets of $\ck(W)$, $Val(W)$ is a
Banach space (see e.g. Lemma 7.0.3 in \cite{alesker-Kent}).

\item \begin{definition}\label{D:homog}
A valuation $\phi\in Val(W)$ is called $i$-homogeneous if
$$\phi(\lam K)=\lam^i\phi(K) \mbox{ for any } \lam>0,\, K\in \ck(W).$$
\end{definition}

Let $Val_i(W)\subset Val(W)$ denote the subspace of $i$-homogeneous valuations; clearly it is a closed linear subspace.
The following result is called McMullen's decomposition theorem.
\begin{theorem}\label{T:mcmullen-decomp}[McMullen \cite{mcmullen-decomp}]
$$Val(W)=\oplus_{i=0}^{\dim W}Val_i(W),$$
where the sum runs over all integers between 0 and $\dim W$.
\end{theorem}

It is easy to see that 0-homogeneous valuations $Val_0(W)$ is a 1-dimensional space spanned by the Euler characteristic $\chi$.\footnote{By definition $\chi$ takes value 1 on each non-empty convex compact set.}
$Val_{\dim W}(W)$ is also 1-dimensional and is spanned by a Lebesgue measure; this result is due to Hadwiger \cite{hadwiger-vol} (see also \cite{alesker-Kent}, Theorem 3.1.1).

\item The group $GL(W)$ of invertible linear transformations acts linearly and continuously on $Val(W)$ as
$$(g\phi)(K)=\phi(g^{-1}K) \mbox{ for any } g\in GL(W), \phi\in Val(W), K\in \ck(W).$$
\begin{definition}\label{D:smooth}
A valuation $\phi\in Val(W)$ is called {\itshape smooth} if the map $GL(W)\to Val(W)$ given by
$$g\mapsto g\phi$$
is $C^\infty$-differentiable.
\end{definition}
Denote the subset of all smooth valuations by $Val^\infty(W)$. A well known representation theoretical result is that $Val^\infty(W)$ is a $GL(W)$-invariant linear subspace dense in $Val(W)$.
Moreover $Val^\infty(W)$ admits a canonical Fr\'echet topology, called the Garding topology, which is stronger than that induced from $Val(W)$. We do not need the precise definition of it, but it is important to know that
the product and the convolutions on $Val^\infty(W)$  discussed below are continuous with respect to Garding topology. Clearly $Val^\infty(W)$ also satisfies the McMullen's decomposition theorem.

 An important example of smooth valuations which will be used below is $K\mapsto vol(K+A)$ where $A\in \ck(W)$ has infinitely smooth boundary and positive Gauss curvature
(i.e. all principal curvatures of the boundary are strictly positive; this class of convex bodies is independent of any Euclidean metric). Here $K+A$ is the Minkowski sum
$$K+A:=\{k+a|k\in K,\, a\in A\}.$$

Furthermore the mixed volume
$K\mapsto V(K [i], A_1,\dots,A_{n-i})$
is a smooth valuation provided $A_1,\dots A_{n-i}\in \ck(W)$ have smooth boundary and positive Gauss curvature. Here and below the notation $K[i]$ means that the body $K$ is repeated $i$ times.

\item Let us discuss the product on valuations. We denote by $\iota_1,\iota_2,\Delta\colon W\to W\times W$ the imbeddings
given by $\iota_1(w)=(w,0),\, \iota_2(w)=(0,w),\, \Delta(w)=(w,w)$.
\begin{theorem}\label{T:product}[\cite{alesker-product}]
(1) There exists a continuous (in the Garding topology) bilinear map called product
$$Val^\infty(W)\times Val^\infty(W)\to Val^\infty(W)$$
which is uniquely characterized by the following property: Let $A,B\in \ck(W)$ have smooth boundary and positive curvature. Consider smooth valuations
$\phi_A(\bullet)=vol_W(\bullet+A),\phi_B(\bullet)=vol_W(\bullet+B)$ where $vol_W$ is a Lebesgue measure on $W$. Then their product is
$$(\phi_A\cdot\phi_B)(K)=vol_{W^2}(\Delta(K)+(A\times B)),$$
where $vol_{W^2}=vol_W\times vol_W$ is the product measure.

(2) Equipped with this product $Val^\infty(W)$ is an associative commutative algebra with the unit $\chi$ (the Euler characteristic).

(3) $Val^\infty(W)$ is a graded algebra with respect to the McMullen's decomposition, i.e. $Val_i^\infty(W)\cdot Val^\infty_j(W)\subset Val_{i+j}^\infty(W).$

(4) (Poincar\`e duality) Consider the product map of smooth valuations of complementary degree of homogeneity
$$Val_i^\infty(W)\times Val_{\dim W-i}^\infty(W)\to Val^\infty_{\dim W}(W)=\CC\cdot vol_W.$$
This is a perfect pairing, i.e. for any $0\ne\phi\in Val_i^\infty(W)$ there exists $\psi \in Val_{\dim W-i}^\infty(W)$ such that
$$\phi\cdot\psi \ne 0.$$
\end{theorem}

\begin{example}\label{EX-mixed-vol-prod}
Let $\phi(\bullet)=V(\bullet[i], A_1,\dots,A_{n-i}),\, \psi(\bullet)=V(\bullet[j], B_1,\dots,B_{n-j})$ where $i+j\leq n$ and $A_p,B_q\in\ck(W)$ have smooth boundary with positive Gauss curvature.

(1) Then
\begin{eqnarray}\label{E:mixed-prod-formula}
(\phi\cdot \psi)(\bullet)={2n\choose n}{i+j\choose i}^{-1}V(\Delta(\bullet)[i+j]; \iota_1A_1,\dots,\iota_1A_{n-i};\iota_2B_1,\dots,\iota_2B_{n-j}),
\end{eqnarray}
where the mixed volume in the right hand side is taken in $W\times W$.

(2) If $i+j=n$ then the expression for the product can be presented in a more explicit form by Proposition 2.2 in \cite{alesker-product}:
\begin{eqnarray}\label{E:product-complement-degree}
(\phi\cdot \psi)(\bullet)={n \choose i}^{-1}V(A_1,\dots,A_{n-i},-B_1,\dots,-B_{i}) vol_W(\bullet),
\end{eqnarray}
where the mixed volume is taken in $W$.
\end{example}

\item Let us discuss convolution on valuations which was introduced by Bernig and Fu.
\begin{theorem}\label{T;convolution}[\cite{bernig-fu-convolution}]
Let us fix a positive Lebesgue measure $vol_W$ on the space $W$.

(1) There exists a continuous (in the Garding topology) bilinear map called convolution\footnote{The convolution does depend on a choice of $vol_W$.}
$$Val^\infty(W)\times Val^\infty(W)\to Val^\infty(W)$$
which is uniquely characterized by the following property: Let $A,B\in \ck(W)$ have smooth boundary and positive Gauss curvature. Consider smooth valuations
$\phi_A(\bullet)=vol_W(\bullet+A),\phi_B(\bullet)=vol_W(\bullet+B)$. Then
$$(\phi_A\ast\phi_B)(\bullet)=vol_W(\bullet+A+B).$$

(2) Equipped with the convolution $Val^\infty(W)$ is an associative commutative algebra with unit $vol_W$.

(3) $Val^\infty_i(W)\ast Val^\infty_j(W)\subset Val^\infty_{i+j-\dim W}(W).$
\end{theorem}

It was shown in \cite{alesker-fourier} that the topological algebra $Val^\infty(W)$ equipped with the product is isomorphic to the topological algebra $Val^\infty(W^*)$ equipped with convolution\footnote{$W^*$ denotes the dual space of $W$.}; moreover an
isomorphism may be chosen to commute with the natural action of the subgroup of $GL(W)$ of volume preserving transformations.

\begin{example}\label{EX-mixed-vol-colv-conv}
Let $\phi(\bullet)=V(\bullet[i], A_1,\dots,A_{n-i}),\, \psi(\bullet)=V(\bullet[j], B_1,\dots,B_{n-j})$ where $i+j\geq n$ and $A_p,B_q\in\ck(W)$ have smooth boundary with positive Gauss curvature.
Then
\begin{eqnarray}\label{EX-convol-mixed}
(\phi\ast\psi)(\bullet)={i+j\choose n}{i+j\choose i}^{-1}V(\bullet [i+j-n];A_1,\dots,A_{n-i};B_1,\dots,B_{n-j}),
\end{eqnarray}
where the mixed volume in the right hand side is in $W$.

\end{example}

\item Let us discuss hard Lefschetz type theorems on valuations. There are two versions: for the product and for the convolution. We denote by $V_i$, $i=0,1,\dots, n$, the $i$th intrinsic volume.
Recall that $V_i(\bullet)$ is proportional with positive coefficient
to the mixed volume with the unit  Euclidean ball $V(\bullet[i],B[n-i])$. The product $V_i\cdot V_j$ is proportional with positive coefficient to $V_{i+j}$, and the convolution $V_i\ast V_j$ is proportional with
positive coefficient to $V_{i+j-n}$.
\begin{theorem}\label{T:hard-Lefschetz}
Let us denote $n=\dim W$. Let $0\leq i< n/2$.

(1)  The linear map
$Val^\infty_i(W)\to Val_{n-i}^\infty(W)$
given by $\phi\mapsto \phi\cdot (V_1)^{n-2i}$ is an isomorphism.

(2) The linear map $Val_{n-i}^\infty(W)\to Val_i^\infty(W)$ given by $\psi\mapsto \psi\ast (V_{n-1})^{\ast(n- 2i)}$ is an isomorphism.

\end{theorem}

Initially the version (2) was proven by the author \cite{alesker-jdg-03} for even valuations only. Version (2) in general was proven by Bernig and Br\"ocker \cite{bernig-brocker}.
Version (1) was proven for even valuations by the author \cite{alesker-hard-lefschetz-product}. Then it was proven in general in \cite{alesker-fourier} by the author; the proof used version (2) by Bernig and Br\"ocker \cite{bernig-brocker}
and a version of the Fourier transform on valuations introduced in \cite{alesker-fourier}.

\item Let us state Kotrbat\'y's results. There will be two equivalent versions formulated either in terms of product or convolution.

\begin{definition}\label{D:primitive}
Let $i\leq n/2$ where $n=\dim W$.

(1) A valuation $\phi\in Val_i^\infty(W)$ is called {\itshape primitive} if $\phi\cdot (V_1)^{n-2i+1}=0.$

(2) A valuation $\psi\in Val_{n-i}^\infty(W)$ is called {\itshape co-primitive} if $\psi\ast (V_{n-1})^{\ast(n- 2i+1)}=0$.
\end{definition}
For $i\leq n/2$ let us define the Hermitian form $Q$ on $Val_i^\infty(W)$ with values in $Val_n(W)=\CC\cdot vol_W$ by
$$Q(\phi):=(-1)^i\phi\cdot\bar\phi\cdot (V_1)^{n-2i}$$
and the Hermitian form $\tilde Q$ on $Val_{n-i}^\infty(W)$ with values in $Val_0(W)=\CC\cdot \chi$ by
$$\tilde Q(\psi)=(-1)^i \psi\ast\bar\psi \ast V_{n-1}^{\ast(n- 2i)}.$$

Kotrbat\'y has formulated the following conjecture which is an analogue of the Hodge-Riemann bilinear relations from K\"ahler geometry.
\begin{conjecture}\label{kotrbaty-conj}[\cite{kotrbaty-hr}]
Let $i\leq n/2$.

(1) Let $\phi\in Val_i^\infty(W)$ be a non-zero primitive even (resp. odd) valuation. Then $Q(\phi)>0$ (resp. $Q(\phi)<0$).

(2) Let $\psi\in Val_{n-i}^\infty(W)$ be a non-zero co-primitive valuation. Then $\tilde Q(\psi)>0$.
\end{conjecture}
Kotrbat\'y has shown in \cite{kotrbaty-hr} that parts (1) and (2) of the conjecture are equivalent. He proved the conjecture in the following special cases.
\begin{theorem}\label{T:kotrbaty}[\cite{kotrbaty-hr}]
Conjecture \ref{kotrbaty-conj} holds for even valuations for any $i$ as in the conjecture, and for odd valuations for $i=0,1$.
\end{theorem}

This theorem will be applied to the main results of this paper for valuations given by mixed volumes. Their connection to the product and the convolution
is given in Examples \ref{EX-mixed-vol-prod} and \ref{EX-mixed-vol-colv-conv}.

\end{paragraphlist}

\section{Proof of the main results.}\label{S:pf-main}
Let $W$ be an $n$-dimensional vector space with a fixed positive Lebesgue measure $vol_W$. We denote the linear maps
\begin{eqnarray*}
\Delta_2\colon W\inj W^2 \mbox{given by } \Delta_2 (w)=(w,w),\\
\Delta_3\colon W\inj W^3 \mbox{given by } \Delta_3 (w)=(w,w,w).
\end{eqnarray*}
As previously we denote the imbeddings $\iota_1,\iota_2\colon W\inj W\times W$ by $\iota_1(w)=(w,0)$, $\iota_2(w)=(0,w)$. In a similar way, by the abuse of notation, we will denote analogous imbeddings
$\iota_1,\iota_2,\iota_3\colon W\inj W^3$ as imbeddings of $W$ into the corresponding copy of $W$ in the triple  product $W^3$,
\begin{lemma}\label{L:triple-product}
Let $A,B,C\subset W$ be convex compact subsets with smooth boundary and positive Gauss curvature. Consider smooth valuations
$$\phi_A(\bullet):=vol_W(\bullet +A),$$
and similarly for $B,C$. Then
$$(\phi_A\cdot\phi_B\cdot \phi_C)(K)= vol_{W^3}(\Delta_3(K)+(A\times B\times C)) \mbox{ for any } K\in\ck(W),$$
where $vol_{W^3}$ is the product measure of $vol_W$ taken 3 times.
\end{lemma}
\begin{remark}
Analogous statement can be proven not only for a triple product of valuations of the above form, but for any number of them.
\end{remark}

{\bf Proof.} First let us show that for any $\psi\in Val^\infty(W)$ and $\phi_C$ as in the lemma one has
\begin{eqnarray}\label{E:product-smooth-00}
(\psi\cdot \phi_C)(K)=\int_W \psi(K\cap(y-C))dy.
\end{eqnarray}
For by continuity and by solution of the McMullen's conjecture \cite{alesker-gafa-01} it suffices to assume that $\psi$ is of the form $\phi_A$. Then by Fubini theorem
\begin{eqnarray*}
(\phi_A\cdot \psi_C)(K)=vol_{W^2}(\Delta_2(K)+(A\times C))=\\
\int_W vol_W\left(\left[\Delta_2(K)+(A\times C)\right]\cap (W\times\{y\})\right) dy=\\
\int_W vol_W\left(\left(\left[\Delta_2(K)+(0\times C)\right]\cap (W\times\{y\}) \right)+A\right)dy=\\
\int_W \phi_A\left(\left[\Delta_2(K)+(0\times C)\right]\cap (W\times\{y\})\right) dy=\\
\int_W \phi_A\left(\{k|\, k\in K \mbox{ and }\exists c\in C \mbox{ s.t. } k+c=y\} \right)dy=\\
\int_W \phi_A(K\cap (y-C))dy
\end{eqnarray*}
which is (\ref{E:product-smooth-00}).

Now let us consider
\begin{eqnarray*}
vol_{W^3}(\Delta_3(K)+(A\times B\times C))=\\
\int_W vol_{W^2}\left(\left(\left[\Delta_3(K)+(0\times 0\times C)\right]\cap (W^2\times \{y\})\right)+(A\times B\times 0\right) dy=\\
\int vol_{W^2}\left(\left\{(k,k)|\, k\in K\mbox{ and }\exists c\in C \mbox{ s.t. } k+c=y\right\}+(A\times B)\right)dy=\\
\int vol_{W^2}\left(\Delta_2(K\cap (y-C))+(A\times B)\right)dy=\\
\int (\phi_A\cdot\phi_B)(K\cap (y-C))dy\overset{(\ref{E:product-smooth-00})}{=}
((\phi_A\cdot\phi_B)\cdot\phi_C)(K).
\end{eqnarray*}
\qed

\begin{lemma}\label{L:technical-computtion}
Let us given two $n-1$-tuples of convex compact sets $A_i,B_i,\, i=1,\dots,n-1,$ in $n$-dimensional space $W$ , and two more convex compact sets $C_1,C_2$. Let us assume that all the bodies have smooth boundary and
positive curvature. Consider valuations
\begin{eqnarray*}
\alp(\bullet)=V(\bullet, A_1,\dots,A_{n-1}),\\
\beta(\bullet)=V(\bullet,B_1,\dots,B_{n-1}),\\
\gamma(\bullet)=V(\bullet[n-2],C_1,C_2).
\end{eqnarray*}
Then their product is equal to
$$\alp\cdot\beta\cdot\gamma=c_n V(\iota_1A_1,\dots,\iota_1A_{n-1},\iota_2B_1,\dots,\iota_2B_{n-1}, -\Delta_2(C_1),-\Delta_2(C_2))\cdot vol_W,$$
where $c_n>0$ is a constant depending on $n$ only.\footnote{$c_n$ can be computed explicitly, but its value is not important in this paper.}
\end{lemma}
{\bf Proof}. Clearly
\begin{eqnarray*}
\alp(\bullet)=\frac{1}{n!}\frac{\pt}{\pt \lam_1}\dots \frac{\pt}{\pt \lam_{n-1}}vol_W(\bullet+\sum_{i=1}^{n-1}\lam_iA_i),\\
\beta(\bullet)=\frac{1}{n!}\frac{\pt}{\pt \mu_1}\dots \frac{\pt}{\pt \mu_{n-1}}vol_W(\bullet+\sum_{j=1}^{n-1}\mu_jB_j),\\
\gamma(\bullet)={n\choose 2}^{-1}\frac{\pt}{\pt\eps_1}\frac{\pt}{\pt\eps_2}vol_W(\bullet+\eps_1C_1+\eps_2C_2),
\end{eqnarray*}
where all the partial derivatives are taken at 0. Hence by Lemma \ref{L:triple-product} we get
\begin{eqnarray}\label{E:abetagamma}
(\alp\cdot\beta\cdot\gamma)(K)=c_n' V(\Delta_3(K)[n];\iota_1A_1,\dots,\iota_1A_{n-1};\iota_2B_1,\dots,\iota_2B_{n-1};\iota_3C_1,\iota_3C_2),
\end{eqnarray}
where $K\in \ck(W)$ is arbitrary, and $c_n'$ is a constant depending on $n$ only. The product is $n$-homogenous valuation in $K$ and hence has to be proportional to $vol_W(K)$
by the Hadwiger theorem \cite{hadwiger-vol}. However this fact can be seen more directly using Theorem 5.3.1 in \cite{schneider-2ed} which will also give the required coefficient of proportionality.

The latter result says that if $E\subset X$ be a $k$-dimensional linear subspace of an $N$-dimensional vector space $X$. Let $\pi\colon X\to X/E$ be the quotient map. Let $L_1,\dots L_{N-k}\in \ck(X)$. Then for any $K\subset E$
$$V(K[k],L_1,\dots,L_{N-k})=c'' vol_E(K) V(\pi L_1,\dots,\pi L_{N-k}),$$
where $vol_E$ is a Lebegue measure on $E$, $V$ in the left (resp. right) hand side is the mixed volume with respect to a Lebesgue measure on $X$ (resp. $X/E$), and $c''>0$ depends on the choice of Lebesgue measures on $E,X,X/E$, its precise value is not important for us.

We apply this result for $E=\Delta_3(W)\subset X=W^3$. Then this theorem and (\ref{E:abetagamma}) imply
\begin{eqnarray}\label{E:abetagamma2}
(\alp\cdot\beta\cdot\gamma)(K)=c_n'' V(\pi\iota_1A_1,\dots \pi \iota_1A_{n-1};\pi\iota_2B_1,\dots \pi \iota_2B_{n-1};\pi\iota_3C_1,\pi\iota_3 C_2)\cdot vol_W(K),
\end{eqnarray}
where $\pi\colon W^3\to W^3/\Delta_3(W)$ is the canonical quotient map. There exists a unique isomorphism of vector spaces
$$Q\colon W^3/\Delta_3(W)\tilde\to W^2,$$
given by $(Q\circ\pi)(w_1,w_2,w_3)=(w_1-w_3,w_2-w_3)$. Applying $Q$ to all convex bodies $A_i,B_i,C_i$ and taking into account that $(Q\pi\iota_1)(A)=\iota_1A,\, (Q\pi\iota_2)(B)=\iota_2B, (Q\pi\iota_3)(C)=-\Delta_2(C)$, one gets
\begin{eqnarray*}
(\alp\cdot\beta\cdot\gamma)(K)=\\
c_n V((Q\pi\iota_1)A_1,\dots (Q\pi \iota_1)A_{n-1};(Q\pi\iota_2)B_1,\dots (Q\pi \iota_2)B_{n-1};(Q\pi\iota_3)C_1,(Q\pi\iota_3) C_2)\cdot vol_W(K)=\\
c_n V(\iota_1A_1,\dots,\iota_1A_{n-1};\iota_2B_1,\dots,\iota_2B_{n-1};-\Delta_2(C_1),-\Delta_2(C_2))\cdot vol_W(K),
\end{eqnarray*}
where the last mixed volume is in $W^2$. \qed

\begin{theorem}\label{T:odd-ineq}
(a) Let $A_1,\dots A_{n-1}\subset\RR^n$ be convex compact subsets of the Euclidean space. Let $B\subset \RR^n$ be the unit Euclidean ball. Then one has inequality for mixed volumes in $\RR^{2n}=\RR^n\times\RR^n$
\begin{eqnarray*}\label{E:odd-ineq}
V(\iota_1A_1,\dots,\iota_1A_{n-1};\iota_2A_1,\dots,\iota_2A_{n-1};\Delta_2(B)[2])\geq\\ V(\iota_1A_1,\dots,\iota_1A_{n-1};-\iota_2A_1,\dots,-\iota_2A_{n-1};\Delta_2(B)[2]).
\end{eqnarray*}

(b) In the special case when $A_1=\dots=A_{n-1}=:A$ with $A$ having non-empty interior and smooth boundary with positive Gauss curvature the equality in the above inequality
is achieved if and only if $A$ has a center of symmetry.
\end{theorem}
{\bf Proof.} We may and will assume that all bodies $A_i$ have smooth boundary and positive curvature. Consider the valuation
$$\phi(\bullet)=V(\bullet, A_1,\dots,A_{n-1})-V(\bullet,-A_1,\dots,-A_{n-1}).$$
$\phi$ is a 1-homogeneous odd valuation, and hence it is primitive. Hence by Kotrbat\'y's theorem \ref{T:kotrbaty}
\begin{eqnarray}\label{E:hr-odd-applic}
\phi^2\cdot V(\bullet[n-2], B[2])\geq 0.
\end{eqnarray}
Computing this product explicitly using Lemma \ref{L:technical-computtion} one gets the result (a).

Note that if $\phi\ne 0$ then there is a strict inequality in (\ref{E:hr-odd-applic}). Hence if there is an equality in the inequality of part (a) of the theorem it follows that $\phi=0$.
In the case $A_1=\dots=A_{n-1}=A$ that is equivalent to the equality of the area measures
$$S_{n-1}(A,\cdot)=S_{n-1}(-A,\cdot).$$
By Theorem 8.1.1 in \cite{schneider-2ed} the latter condition is satisfied if and only if $A$ and $-A$ are translates of each other. This implies part (b) of the theorem.
\qed

\begin{theorem}\label{T:even-ineq}
Let $n\geq 2$. Let $A_1,\dots A_{n-1}\subset\RR^n$ be convex compact subsets of  Euclidean space. Let $B\subset \RR^n$ be the unit Euclidean ball. 

(a) One has inequality
for mixed volumes in $\RR^{2n}=\RR^n\times\RR^n$ in the left hand side and in $\RR^n$ in the right hand side
\begin{eqnarray*}
V(\iota_1A_1,\dots,\iota_1A_{n-1};\iota_2A_1,\dots,\iota_2A_{n-1};\Delta_2(B)[2])+\\
+V(\iota_1A_1,\dots,\iota_1A_{n-1};-\iota_2A_1,\dots,-\iota_2A_{n-1};\Delta_2(B)[2])\leq \gamma_n V(A_1,\dots,A_{n-1},B)^2,
\end{eqnarray*}
where $\gamma_n$ is a constant depending on $n$ only and uniquely characterized by the property that if $A_1=\dots=A_{n-1}$ are equal to $B$ then in the above inequality there is an equality.

(b) In the special case when $A_1=\dots=A_{n-1}=:A$ with $A$ being a centrally symmetric convex body with non-empty interior and smooth boundary with positive Gauss curvature, the equality in the above inequality
is achieved if and only if $A$ is a Euclidean ball of arbitrary radius.
\end{theorem}
{\bf Proof.} We may and will assume that all convex bodies $A_i$ have smooth boundary with positive curvature. Consider the valuation
$$\phi(\bullet):=V(\bullet,A_1,\dots,A_{n-1})+V(\bullet,-A_1,\dots,-A_{n-1}),$$
Clearly $\phi$ is an even 1-homogeneous valuation.
We are going to show that there exists $\lam\in \RR$ such that the valuation
$$\psi(\bullet):=\phi(\bullet)-\lam V(\bullet,B[n-1])$$
is primitive. That means that
$0=\psi\cdot V(\bullet[n-1],B),$ or equivalently
\begin{eqnarray}\label{E000}
\phi\cdot V(\bullet[n-1],B)=\lam V(\bullet,B[n-1])\cdot V(\bullet[n-1],B).
\end{eqnarray}

By (\ref{E:product-complement-degree}) this is equivalent to
\begin{eqnarray*}
2V(B,A_1,\dots,A_{n-1})=\lam vol(B).
\end{eqnarray*}
Thus
\begin{eqnarray}\label{E:lambda}
\lam=2\frac{V(B,A_1,\dots,A_{n-1})}{vol(B)}.
\end{eqnarray}

For this particular $\lam$ let us use Kotrbat\'y's theorem \ref{T:kotrbaty}:
\begin{eqnarray*}
0\geq \psi^2\cdot V(\bullet[n-2],B[2])=\phi^2V(\bullet[n-2],B[2])+\lam^2 V(\bullet,B[n-1])^2\cdot V(\bullet[n-2],B[2])-\\
-2\lam \phi V(\bullet,B[n-1])\cdot V(\bullet[n-2],B[2])
\end{eqnarray*}
By (\ref{E000}) the last summand is equal to $-2\lam^2V(\bullet,B[n-1])^2\cdot V(\bullet[n-2],B[2])$. Hence we obtain
\begin{eqnarray}\label{E:inequality-even-etc}
0\geq \phi^2V(\bullet[n-2],B[2])-\lam^2 V(\bullet,B[n-1])^2\cdot V(\bullet[n-2],B[2])
\end{eqnarray}

By Lemma (\ref{L:technical-computtion})
\begin{eqnarray*}
 \phi^2V(\bullet[n-2],B[2])=c_nvol\cdot (V(\iota_1A_1,\dots,\iota_1 A_{n-1},\iota_2A_1,\dots,\iota_2 A_{n-1},\Delta_2(B)[2])+\\
 V(-\iota_1A_1,\dots,-\iota_1 A_{n-1},-\iota_2A_1,\dots,-\iota_2 A_{n-1},\Delta_2(B)[2])+\\2V(\iota_1A_1,\dots,\iota_1 A_{n-1},-\iota_2A_1,\dots,-\iota_2 A_{n-1},\Delta_2(B)[2]))=\\
 2c_nvol\cdot(V(\iota_1A_1,\dots,\iota_1 A_{n-1},\iota_2A_1,\dots,\iota_2 A_{n-1},\Delta_2(B)[2])+\\V(\iota_1A_1,\dots,\iota_1 A_{n-1},-\iota_2A_1,\dots,-\iota_2 A_{n-1},\Delta_2(B)[2])).
\end{eqnarray*}

Similarly
\begin{eqnarray*}
V(\bullet,B[n-1])^2\cdot V(\bullet[n-2],B[2])=c_n vol\cdot V(\iota_1B[n-1],\iota_2B[n-1],\Delta_2[2]).
\end{eqnarray*}

Substituting the last two equalities into (\ref{E:inequality-even-etc}) we get
\begin{eqnarray*}
 V(\iota_1A_1,\dots,\iota_1 A_{n-1},\iota_2A_1,\dots,\iota_2 A_{n-1},\Delta_2(B)[2])+\\V(\iota_1A_1,\dots,\iota_1 A_{n-1},-\iota_2A_1,\dots,-\iota_2 A_{n-1},\Delta_2(B)[2])\leq \\
\frac{\lam}{2}^2V(\iota_1B[n-1],\iota_2B[n-1],\Delta_2(B)[2])\overset{(\ref{E:lambda})}{=}\\
(V(B,A_1,\dots,A_{n-1}))^2\underset{\gamma_n}{\underbrace{\frac{2V(\iota_1B[n-1],\iota_2B[n-1],\Delta_2(B)[2])}{vol(B)^2}}}.
\end{eqnarray*}
Thus part (a) of the theorem is proven. 

To prove part (b) assume that $A=-A$ and there is an equality in the inequality of part (a). This is equivalent to $\psi=0$, or 
$$V(\bullet, A[n-1])=\frac{\lam}{2}V(\bullet,B[n-1]).$$
In other words $S_{n-1}(A,\cdot)=\frac{\lam}{2}S_{n-1}(B,\cdot)$. By Theorem 8.1.1 in \cite{schneider-2ed} $A=\frac{\lam}{2}B$.
\qed

\hfill

If we combine Theorems \ref{T:odd-ineq} and \ref{T:even-ineq} we get
\begin{corollary}\label{Cor-ineq}
Let $n\geq 2$. Let $A_1,\dots A_{n-1}\subset\RR^n$ be convex compact subsets of  Euclidean space. Let $B\subset \RR^n$ be the unit Euclidean ball. Then one has
\begin{eqnarray*}
V(\iota_1A_1,\dots,\iota_1A_{n-1};-\iota_2A_1,\dots,-\iota_2A_{n-1};\Delta_2(B)[2])\leq \gamma_n' \left(V(A_1,\dots,A_{n-1},B)\right)^2,
\end{eqnarray*}
where $\gamma_n'$ is a constant depending on $n$ only and uniquely characterized by the property that if $A_1=\dots=A_{n-1}=B$  then in the above inequality there is an equality.\footnote{In fact $\gamma_n'=\gamma_n/2$
where $\gamma_n$ is the constant from Theorem \ref{T:even-ineq}.}
\end{corollary}

\begin{theorem}\label{T:mixed-area}
Let $n\geq 4$. Let $K_1,K_2,L_1,L_2\subset \RR^n$ be centrally symmetric convex compact sets with smooth boundary with positive Gauss curvature. Let $B\subset \RR^n$ denote the unit Euclidean ball.
Let us assume the following equality of the mixed area measures:
\begin{eqnarray}\label{E:mixed-meas-cond}
S(K_1,K_2,B[n-3],\cdot)=S(L_1,L_2,B[n-3],\cdot).
\end{eqnarray}
Then one has
$$V(K_1[2],K_2[2],B[n-4])+V(L_1[2],L_2[2],B[n-4])\geq 2V(K_1,K_2,L_1,L_2,B[n-4]).$$
\end{theorem}
{\bf Proof.} Let us consider $Val^\infty(\RR^n)$ as algebra equipped with the Bernig-Fu convolution. Consider the even smooth $n-2$ homogeneous valuation
$$\phi(\bullet)=V(\bullet[n-2],K_1,K_2)-V(\bullet[n-2],L_1,L_2).$$
Let us show that it is primitive in the sense of convolution, i.e. $\phi\ast V_{n-1}^{n-3}=0$, or equivalently
$\phi\ast V(\bullet[3],B[n-3])=0$. By Example \ref{EX-mixed-vol-colv-conv} this is equivalent to
$$V(\bullet,K_1,K_2,B[n-3])-V(\bullet,L_1,L_2,B[n-3])=0.$$
This equality is equivalent to the assumption (\ref{E:mixed-meas-cond}).

By the Kotrbat\'y's Theorem \ref{T:kotrbaty} and by Example \ref{EX-mixed-vol-colv-conv} one has for a positive constant $\kappa_n>0$
\begin{eqnarray*}
0\leq \phi^2\ast V_{4}=\kappa_n vol\cdot (V(K_1[2],K_2[2],B[n-4])+\\
+V(L_1[2],L_2[2],B[n-4])-2V(K_1,K_2,L_1,L_2,B[n-4])).
\end{eqnarray*}
This is exactly the statement of the theorem. \qed

\end{document}